\documentclass{article}
\usepackage{amsmath,amssymb,amsfonts,amsthm}

\usepackage{graphicx}

\newtheorem{theorem}{Theorem}
\newtheorem{lemma}{Lemma}
\newtheorem{corollary}{Corollary}
\begin{document}

\title{On orthogonal curvilinear coordinate systems in constant curvature spaces}
%\thanks{The work is %supported by RFBR (grant %09--01--00598 (I.P.Rybnikov)
%and grant 10--01--91056 %(D.A.Berdinsky)) and Council %of president of
%the Russian Federation for %research grants and support %of leading
%scientific school (grant %7256.2010.1).}
\author{Berdinsky~D.
\thanks{email: berdinsky@gmail.com.}
%Department of Mechanics and %Mathematics, Novosibirsk %State University,
%630090 Novosibirsk, Russia; 
\and Rybnikov~I.
\thanks{email: ivan.p.rybnikov@gmail.com.
%Department of Mechanics and %Mathematics, Novosibirsk %State University,
%630090 Novosibirsk, Russia; 
} }

\date{}
\maketitle {\small
\begin{quote}
\noindent{ {\bf Abstract}.  We describe a method for constructing $n$--orthogonal  coordinate systems in constant
curvature spaces. The construction proposed is a modification of Krichever's method for
producing orthogonal curvilinear coordinate systems in the $n$--dimensional Euclidean space. To demonstrate
how this method works, we construct  examples of orthogonal coordinate systems on the 
two--dimensional sphere and the hyperbolic plane, in the case when the spectral curve is reducible and all irreducible components are isomorphic to a complex projective line.
 } 

\noindent{\bf Keywords}: orthogonal coordinate systems, spaces of constant curvature, Baker--Akhiezer function 
  \end{quote}

\section{Introduction}

This paper presents a method for constructing orthogonal  
coordinate systems in the spaces of constant curvature  
$K \neq 0$ in terms of $n$--point Baker--Akhiezer functions.

We briefly recall some results regarding the classical problem of constructing
orthogonal curvilinear coordinate systems in $\mathbb R ^n$ $(K = 0)$.
Let us consider an orthogonal coordinate system $u_1,\dots,u_n$, and let the metric have
the form:
$$
ds^2 = H_1 ^2 (du^1 ) ^2 + \dots + H_n ^2 (du^n ) ^2,
$$
where $H_i (u^1, \dots ,u^n)$ are the 
Lam\'{e} coefficients. Then the metric is
flat if and only if:
\begin{equation}
\partial _k \beta_{ij} = \beta_{ik} \beta_{kj},\quad i \neq j \neq k,
\end{equation}
\begin{equation}
\partial _{i} \beta_{ij} + \partial_j \beta_{ji} + \sum_{m \neq i,j} \beta_{mi} \beta_{mj} = 0,
\end{equation}
where $\beta_{ij}$ are the rotation coefficients, $ \beta_{ij}
=\frac{\partial _i H_j}{H_i}$, $i \neq j$.

%The solution of  (1), (2) %using the methods of soliton
%theory  was proposed by %Zakharov [1].
Zakharov proposed a solution  
to the system (1), (2) using the methods of soliton theory [1].  
%The system (1) has form of
%the abstract $n$--wave problem. 
%The dressing procedure was %applied
%to the system (1) and the %proper differential %reduction was found in [1]
%for solving (1), (2). 
In [1] the dressing procedure 
was applied to the 
system (1) and a proper
differential reduction was found
for solving the entire system (1), (2).  
%The procedure of
%integration  by the inverse %scattering technique was %applied by
%Zakharov in~[2] to the %system describing an %orthogonal coordinate
%system in the spaces of %diagonal curvature (the %spaces of constant
%curvature are the spaces of %diagonal curvature).
In [2] Zakharov applied 
the procedure of integration by the inverse scattering method
to the system 
describing  
orthogonal 
coordinate
systems in spaces of diagonal curvature 
(spaces of constant
curvature are spaces of diagonal curvature).

%The finite--gap integration %method was applied by %Krichever [3] to
%the problem of constructing  %orthogonal curvilinear %coordinate
%systems in $\mathbb R ^n$. 

In [3] Krichever applied the finite--gap integration method for constructing  orthogonal curvilinear coordinate systems in $\mathbb R ^n$. 
In Krichever's method the coordinate
functions $x^j (u_1,\dots,u_n)$, $j=1,\dots , n$ are expressed 
explicitly as: %in an explicit way:
$$x^j (u_1,\dots,u_n) = \psi (u_1,\dots,u_n,Q_j),$$ where $Q_j$ are
 points on a Riemann surface $\Gamma$,
$\psi(u_1,\dots,u_n,z)$ is a $n$--point Baker--Akhiezer function and $z \in \Gamma$.

%The relation of flat %diagonal metrics with %integrable systems of
%hydrodynamical type [4] was %revealed by Tsarev [5] in %1984 which
%revived interest to the %classical problem of %describing
%orthogonal curvilinear %coordinate systems in flat %spaces and its
%applications to mathematical %physics.

A relation between flat diagonal metrics and integrable systems of
hydrodynamic type [4] was found by Tsarev [5]. %in 1984
This 
revived the interest in the classical problem of describing
orthogonal curvilinear coordinate systems in flat spaces and its
applications in mathematical physics.
%By using the flat diagonal %Darboux--Egorov metrics, it %is possible to
%construct solutions of %Witten~--- Dijkgraaf~--- %Verlinde~---
%Verlinde associativity %equations~[3,\,6,\,7].
The flat diagonal Darboux--Egorov metrics
enable to
construct solutions of the Witten--Dijkgraaf--Verlinde--Verlinde associativity equations~[3,\,6,\,7].

%It is easy to verify that %the metric of constant %curvature $K$ corresponds
%to a solution of the %following system:

It can be verified that the metric of constant curvature $K$ corresponds
to a solution of the following system:
\begin{equation}
\partial _k \beta_{ij} = \beta_{ik} \beta_{kj}, \quad i \neq j \neq k,
\end{equation}
\begin{equation}
\partial _{i} \beta_{ij} + \partial_j \beta_{ji}
+ \sum_{m \neq i,j} \beta_{mi} \beta_{mj} = - K H_i H_j, 
\end{equation}
where, as before, $\beta_{ij} = \frac{\partial _i H_j}{H_i}$, $i
\neq j$. The systems (1) and (3) coincide.

%It needs to be emphasized 
We note that
the metrics of constant curvature %emerge
appear in the description of non--local Hamiltonian operators of hydrodynamic
type (see [8], where the condition for a metric to be of constant
curvature is necessary for the Poisson brackets to be skew--symmetric and
to satisfy the Jacobi identity). The non--local Poisson brackets
of hydrodynamic type, generated by 
the metrics of constant curvature 
(Mokhov--Ferapontov brackets), play a significant role in
the theory of systems of hydrodynamic type. The problem of describing
compatible Mokhov--Ferapontov brackets is equivalent to the
problem of describing pencils of metrics of constant curvature. For the latter it
is sufficient to classify pairs of diagonal metrics of constant curvature that have special form, 
that was done by Mokhov in [9].
%This problem was solved in %using the inverse scattering %method [9].
Moreover, Mokhov proved
%earlier 
the integrability of the equations describing pencils of flat
metrics (compatible Dubrovin--Novikov brackets) that correspond to
some essential reductions of the system (1), (2).

Furthermore, in [10] the Lax pairs with a spectral parameter are constructed for %considerably more 
general classes of orthogonal coordinate systems in spaces of constant curvature  
(described by %more 
general integrable reductions of 
the system (3), (4)
and related to the pairs of compatible Poisson brackets of hydrodynamic type, where one of the brackets is a Mokhov--Ferapontov bracket and
another one is an arbitrary non--local Poisson bracket of hydrodynamic type).
 Also, in [11] it is shown that  
 the system (3), (4) %are
 %appears to be  
 is the  compatibility condition for some linear system without a spectral
 parameter. Note that a spectral parameter appears in the integrable reductions
 of the system (3), (4), associated with the compatible non--local Poisson brackets of hydrodynamic type.
 
 In this paper we %specify 
 find the spectral data %such that 
 for which the coordinate
 functions, represented in terms of the corresponding Baker--Akhiezer
 functions, describe orthogonal 
 coordinate systems in spaces of constant curvature.
 We also %give particular 
 show examples of 
 solutions to 
 the system (3), (4). % in an explicit way
  We follow the approach used by 
 Mironov and Taimanov in [7,\,12].

 The rest of the paper is organized as
 follows. 
 In \S\,2 we remind the definition of a $n$--point Baker--Akhiezer
 function. In~\S\,3 we modify Krichever's method [3] to obtain
 orthogonal coordinate systems in $S^n$ and~$H^n$ in terms of
 Baker--Akhiezer functions. In \S\,4 we show examples of
 orthogonal coordinate systems in $S^2$ and $H^2$ for the case of a singular
 reducible spectral curve such that each irreducible component is
 isomorphic to $\Bbb {CP}^1$. In this case  the Baker--Akhiezer functions may be expressed in terms of elementary functions (%whereas 
 in the case
 of a smooth spectral curve they are 
 expressed in terms of the theta functions of the Jacobi variety).

% In~\S\,4 we follow the methods of %[7,\,12].

\section{Multipoint Baker--Akhiezer functions}

In this section we briefly recall the definition of a $n$--point
Baker--Akhiezer function. 
The notion of a $n$--point Baker--Akhiezer
function (for $n=2$) 
%appeared firstly in 
was originally introduced in [13].

Let $\Gamma$ be a Riemann surface of genus $g$. Let
$P,\gamma$ and $R$ be divisors on $\Gamma$:
$$
P=P_1 + \dots + P_n, \quad \gamma = \gamma_1 + \dots + \gamma_{g+l},
\quad R = r + r_1 + \dots + r_l.
$$
Let $k_i ^{-1}$, $i=1, \dots, n$ be 
some local parameters in
neighborhoods of  $P_1,\dots,P_n$.

The {\it $n$--point Baker--Akhiezer function} $\psi ( u_1 , \dots ,
u_n, S) $, where $u_1,\dots, u_n$ are real variables,
$S \in \Gamma$, corresponding to the spectral data $ \{ \Gamma,
P_1,\dots,P_n, k_1,\dots,k_n, \gamma, R \},$ is a function with the
following properties:

\smallskip 

1)~in neighborhoods of $P_1, \dots ,P_n$ the function $\psi$ has
essential singularities of the form:
$$
\psi = e^{k_j u_j} \left( f_j (x,y) + \frac{g_j (x,y)}{k_j } + \dots
\right), \quad j = 1, \dots, n;
$$

2)~the function $\psi$ is meromorphic on $\Gamma \backslash
\{P_1,\dots ,P_n\}$ with simple poles on $\gamma$;

3)~$\psi (u_1,\dots,u_n,r)=h$ and $\psi (u_1,\dots, u_n, r_i) = 0 $,
where $h$ is a nonzero constant (the constants of normalization
for $\psi(u_1,\dots, u_n,r_i), i=1,\dots,l$ 
% are
can be arbitrary, in general, 
but we need
them to be zero for the further proof of Theorem \ref{main_thm} in \S\,3).

\smallskip

%It is possible to express a %Baker--Akhiezer function in terms of a %theta function of $\Gamma$. 

A Baker--Akhiezer function can be 
expressed in terms of a theta function 
for $\Gamma$.   
Let $a_1,\dots,a_g$, $b_1,\dots,b_g$
be a basis of cycles with the following intersection indices:

$$
a_i \circ a_j = b_i \circ b_j = 0,\quad a_i \circ b_j = \delta_{ij},
\quad i,j= 1,\dots, g.
$$
Let $\omega_1, \dots, \omega_g$ be the basis of holomorphic
differentials normalized by the conditions  $\int\limits_{a_j} \omega_i
= \delta_{ij}$. The matrix of $b$--periods $B_{ij}= \int\limits_{b_i}
\omega_j$ is symmetric and its imaginary part is positive definite.

The theta function is defined by the absolutely converging series:
$$
\theta (z) = \sum_{m \in \Bbb {Z}^g} e^{\pi i (Bm,m) + 2\pi i (m,
z)} , \quad z= (z_1, \dots, z_g) \in \Bbb {C} ^g.
$$
For the theta function we have: 
$$
\theta (z+ m) = \theta (z), \quad \theta(z + Bm) = \exp (-\pi i
(Bm,m) - 2 \pi i (m,z)) \theta (z), \quad m \in \Bbb {Z}^g.
$$
Let  $X = \Bbb {C}^g / \{\Bbb {Z}^g + B \Bbb {Z} ^g\}$ be the Jacobi
variety of  $\Gamma$. Let $A: \Gamma \rightarrow X$ be the
Abel map defined as follows:
$$
A(S) = \left( \int\limits_{q_0}^S \omega_1,\dots,
\int\limits_{q_0}^S \omega_g \right), \quad S \in \Gamma,
$$
for some fixed point $q_0$. According to the Riemann theorem, for points
$\gamma_1, \dots ,\gamma_g$ in general position, the function
$\theta(z + A(S))$ has exactly $g$ zeros $\gamma_1,\dots,\gamma_g$
on $\Gamma$, where $z = K_{\Gamma} - A(\gamma_1) - \dots -
A(\gamma_g)$ and $K_\Gamma$ is the vector of Riemann constants.

Let $\Omega ^ i$, $i=1,\dots, n$ be a meromorphic  differential
with pole at  $P_i$ of the form $d \Omega^i = d\big(k_i + O
\big(k_i ^{-1}\big)\big)$ normalized by the conditions
$\int\limits_{a_j} \Omega ^i = 0$, $j=1, \dots,g$. We put $U^i = \big(
\int\limits_{b_1} \Omega ^i , \dots , \int\limits_{b_g}\Omega^ i
\big)$, $i=1,\dots, n$. Let $\widetilde{\psi}$ be a function of the
following form:

\begin{equation}
\begin{split}
 \widetilde{\psi} (u^1,\dots,u^n,S) = \frac{\theta (A(S) +
u^1 U^1 + \dots + u^n U^n +z)}{\theta (A (S) +z)}
\\
\times \exp \left( 2\pi i u_1 \int\limits_{q_0}^S \Omega ^1 + \dots
+ 2 \pi i u^n \int\limits_{q_0}^{S} \Omega^n \right).
\end{split}
\end{equation}
For $l=0$ the Baker--Akhiezer 
function equals:
$$
\psi (u^1,\dots,u^n,S) = f(u^1,\dots,u^n)
\widetilde{\psi}(u^1,\dots,u^n,S)
$$
and $f$ is defined by the 
%condition of normalization 
normalization condition 
$\psi (u^1,\dots,u^n,r) = h$. For $l > 0$, the Baker--Akhiezer function may be %represented 
expressed as 
$ \psi = f \widetilde{\psi} + f_1
\widetilde{\psi}_1 + \dots + f_l \widetilde{\psi}_l$, where
 $\widetilde{\psi}_j, j=1,\dots,l$ are %constructed 
 defined by the divisors $\gamma_1
+ \dots + \gamma_{g-1} +r_j$ analogous to $\widetilde{\psi}$ and
the functions $f$, $f_j$ are defined by 
the normalization conditions. 
%from the conditions of normalization.

\section{Orthogonal %curvilinear
	 coordinate systems for $S^n$ and $H^n$}
\label{modified_krichever_sec}

In this section we introduce 
a modified Krichever's method
for constructing orthogonal  
coordinate systems for the $n$--dimensional  sphere $S^n$ and the hyperbolic space $H^n$.    

We first shortly recall Krichever's method [3]
for constructing orthogonal 
coordinate curvilinear systems in $\mathbb{R}^n$. Let
$\Gamma$ be a smooth algebraic curve of genus $g$ with three
divisors:

$$
P = P_1 + \dots + P_n, \quad D = \gamma_1 + \dots +
\gamma_{g+l-1},\quad R = r_1 + \dots + r_l.
$$
Let us consider a divisor: $$Q=Q_1 + \dots + Q_n,$$ 
where $Q_i \in
\Gamma \setminus \{ P \bigcup D \bigcup R \}$. If we require 
that the spectral curve $\Gamma$ 
and the sets $P$, $D$, $R$ 
%$P\bigcup D \bigcup R$ 
and $Q$ satisfy 
some additional 
conditions, then the functions: 
$$ x^j (u^1,\dots,u^n) = \psi
(u^1,\dots,u^n,Q_j)$$
generate orthogonal coordinates 
in $\Bbb R ^n$. 
These conditions  for 
$\Gamma, P, D, R$ and $Q$ 
are as follows:

\smallskip

 1) %Let   $\Gamma$ admits 
 There exists a holomorphic involution
 $\sigma:\Gamma \rightarrow \Gamma$  such that $P_i$, $Q_j$, $i,j=1,\dots,n$
 are statioinary points of $\sigma$ and $\sigma (k_i ^ {-1}) =- k_i
 ^{-1}$ in some neighborhood of point $P_i$;

 2) There exists a meromorphic differential $\Omega$ on $\Gamma$ with
 the following divisors of zeroes
 and poles:
 $$
  (\Omega)_0 = D + \sigma D + P, (\Omega)_{\infty} = R + 
  \sigma R +  Q
 $$
 %And, moreover, 
 such that 
 $\operatorname{Res}_{Q_1} \Omega = \dots =
 \operatorname{Res}_{Q_n} 
 \Omega = \mu_0 ^2 > 0$;

 3) There exists an antiholomorphic involution $\tau : \Gamma
\rightarrow  \Gamma$ such that $P_i$ and $Q_j$, $i,j=1,\dots,n$
are stationary points of $\tau$, $\tau (k_i ^{-1}) = \overline{k_i ^{-1}}$ in some
neighborhood of $P_i$, $i=1,\dots,n$, 
%In addition, 
$\tau R = R$ and
$\tau D = D$.

\smallskip 

The modified Krichever's method is as follows. Let 
$\Gamma$  be a Riemannian surface of genus $g$ with given
divisors on it:
$$
P = P_1 + \dots + P_n,\quad D = \gamma_1 + \dots + \gamma_{g+l},
\quad R = r + r_1 + \dots + r_l.
$$
Let us take $n+1$ points $Q_1,\dots,Q_{n+1}$ on  $\Gamma$. Assume that 
there exists a holomorphic involution $\sigma : \Gamma
\rightarrow \Gamma$ with stationary points $Q_1, \dots, Q_{n+1}, P_1, \dots ,P_n$ and $r$, 
and in some neighborhood of
 $P_i$ the involution $\sigma$  acts locally as $\sigma (k_i) = -k_i $.
Then %it is easy to verify 
we have the following theorem.
\begin{theorem}
\label{main_thm}
Assume that there exists a 1--form $\Omega$ 
such that
$$
(\Omega)= D + \sigma D + P - Q_1 - \dots -Q_{n+1} - R -\sigma r_1 -
\dots -\sigma r_l.
$$
We use the following notation:
 $A_i = \operatorname{Res}_{Q_i}
\Omega$, $B = \operatorname{Res}_{r} \Omega$ and $C_i$ is defined
from the expansion $\big( -C_i k_i ^{-1} + \dots \big) dk_i ^{-1}$  of
the form $\Omega$ in a neighborhood of  $P_i$, $i=1,
\dots,u_n$. Let $h = \psi (r)$  be the constant of normalization of
the Baker--Akhiezer function $\psi$ at $r$, and let $f_j
(u_1,\dots,u_n)$, $j=1,\dots,n$  be the first coefficient of the expansion
$e^{-u_j k_j} \psi =\big( f_j + \frac{g_j}{k_j}+ \dots \big)$ in a
neighborhood of  $P_j$. Then:
%the following identities hold
\begin{equation}
\label{thm_eq1}
\sum_{k=1} ^{n+1} \psi(Q_k)^2 A_k + h^2 B = 0, 
\end{equation}
\begin{equation}
\label{thm_eq2}
\sum_{k=1} ^{n+1} \psi_{u_i} (Q_k) \psi_{u_j} (Q_k) A_k = 0, \quad
i,j =1,\dots,n,\ i \neq j, 
\end{equation}
\begin{equation}
\sum_{k=1} ^{n+1} \psi_{u_i} ^2 (Q^k) A_k = C_i f_i ^2, \quad
i=1,\dots, n. 
\end{equation}
\end{theorem}
\noindent {\it Proof:} Let us consider the form $\psi (S) \psi(\sigma S)
\Omega$, $S \in \Gamma$. 
%It is easy to see that 
It can be seen that 
this form has poles
only at the points $Q_1,\dots,Q_{n+1}$ and $r$. Since the sum of residues at
these points is %should be 
equal to zero, we get (6).

Now let us consider the form $\psi_{u_i} (S)\psi_{u_j} (\sigma S) \Omega$.
If $i \neq j$, then this form has poles at the points $Q_1,\dots, Q_{n+1}$. Since the sum of residues at these points %should be 
is equal to zero, we get (7). 
If $i=j$, then the poles are located
only at the points $Q_1, \dots, Q_n,P_i$, so we get (8). Theorem is proved. \qed 

\smallskip 

We note that identities analogous to \eqref{thm_eq1} and \eqref{thm_eq2}
already appeared in [14,\,15]. 
The functions $\psi (u_1,\dots,u_n, Q_i )$, $i=1,\dots,n+1$ %depending on 
that depend on $n$ real variables $u_1,\dots,u_n$ are complex--valued.
In order to obtain real--valued functions, we impose the conditions
%of the following lemma
stated in the following lemma. 

\begin{lemma} 
\label{lemma_reality}	
Assume that there exists an antiholomorphic
involution $\tau$ such that the points of divisors $P,D,R$ and the points
$Q_1, \dots, Q_{n+1}$ are stationary with respect to $\tau$. Moreover, in
a neighborhood of $P_i$ the involution $\tau$ has the form $\tau (k_i) =
\overline{k_i}$. Then the functions $\psi (u_1,\dots,u_n, Q_i )$,
$i=1,\dots,n+1$ are real--valued.
\end{lemma}
\noindent {\it Proof:} 
By the uniqueness of a Baker--Akhiezer function, we obtain: $$\psi(u_1,\dots,u_n,S) =
\overline{\psi(u_1,\dots,u_n,\tau S)},$$ which implies that
the functions $\psi(u_1,\dots,u_n,Q_i)$ are real--valued. \qed  %Lemma is proved. 

\smallskip 

%Note that a Baker~--- Akhiezer function of the %type specified in
%\S\,2 and the formulas analogous to (6), (7)  %already appeared in
%[14,\,15], where some surfaces with diagonal %metrics are constructed. 
% topologically different from the sphere

%It is easy now to state the following %corollaries.

Theorem \ref{main_thm} and Lemma \ref{lemma_reality} imply the following corollaries. 

\begin{corollary}
\label{corollary_sphere}
Assume that all $A_i$ are real, 
$A_i > 0$ for all $i=1,\dots, n+1$ and 
$h^2 B =-1$. Then the functions:
$$
x^i (u_1,\dots,u_n) = \sqrt{A_i} \psi \left(u_1,\dots,u_n,Q_i
\right),\quad i=1,\dots,n+1,
$$
give an orthogonal coordinate system in $S^n$. Moreover, the Lam\'{e} coefficients: 
$$H_i =\sqrt{f_i ^2 C_i}$$ 
satisfy the equations  
{\rm (3) and (4)} for $K=1$.
\end{corollary}
\noindent {\it Proof:} The identities 
(6) and (7) in Theorem 
\ref{main_thm}
imply that the functions 
$x^i (u_1,\dots,u_n)$, $i=1,\dots,n+1$
satisfy the  equations:
$$
(x^1) ^2 + \dots + (x^{n+1}) ^2 = 1, \quad x^1 _{u_i} x^1 _{u_j} +
\dots + x^{n+1} _{u_i}x^{n+1} _{u_j} = 0, \quad i \neq j.
$$
%The fact that  $x^i (u_1,\dots,u_n)$, %$i=1,\dots,n+1$ are
%real--valued  follows easily from Lemma. 
It follows from Lemma \ref{lemma_reality}
that  $x^i (u_1,\dots,u_n)$, $i=1,\dots,n+1$
are real--valued. 
%According to (8), we get that 
%$ds^2 = C_1 f_1 ^2 du_1 ^2 + \dots + C_n f_n %^2 du_n ^2$ is a metric of constant curvature %$K=1$. 
%Thus, the Lame coefficients $H_i =
%\sqrt{C_i f_i ^2}$ satisfy (3) and (4). 
Finally, by (8) we obtain that: 
$$ds^2 = C_1 f_1 ^2 du_1 ^2 + \dots + C_n f_n ^2 du_n ^2$$ 
is a metric of constant curvature $K=1$. So the Lam\'{e} coefficients 
$H_i = \sqrt{C_i f_i ^2}$
satisfy the equations (3) and (4). \qed

\smallskip 

%Treating orthogonal coordinate systems in the %hyperbolic space we will
%use the standard model on the hyperboloid, %e.g. imaging $H^n$ as
%the connected component of the hyperboloid %$(x^1) ^2 - (x^2) ^2 - \dots
%- (x^{n+1}) ^2 = 1$ embedded in the Minkowski %space $\Bbb R ^{1,n}$.

As a model of the hyperbolic space $H^n$ we 
use a sheet of the 
two--sheet hyperboloid:  
$$(x^1) ^2 - (x^2) ^2 - \dots
 -(x^{n+1}) ^2 = 1$$ 
embedded in 
the Minkowski space $\Bbb R ^{1,n}$.

\begin{corollary} 
\label{corollary_hyp_space}
Assume that all $A_i$ are real, $A_i < 0$,
for $i=2,\dots,n+1$, $A_1 > 0$ and $h^2 B = -1$. Then the functions:
$$x^i (u^1,\dots,u^n) = \sqrt{|A_i |} \psi ( u^1,\dots,u^n,Q_i ),
\quad i=1,\dots,n+1,
$$
give an orthogonal coordinate system in $H^n$. Moreover, the Lam\'{e}
coefficients: 
$$H_i = \sqrt{-f_i ^2 C_i}$$ 
satisfy  the equations {\rm (3) and (4)} for
$K=-1$.
\end{corollary}
\noindent {\it Proof:} 
The identities (6) and (7) 
imply that the functions $x^i (u_1,\dots,u_n)$,
$i=1,\dots,n+1$ satisfy the equations:
$$
(x^1) ^2 - (x^2)^2 - \dots - (x^{n+1}) ^2 = 1, \quad x^1 _{u_i} x^1
_{u_j} - x^2 _{u_i} x^2 _{u_j} - \dots -x^{n+1} _{u_i} x^{n+1}
_{u_j} = 0, \quad i \neq j.
$$
By Lemma \ref{lemma_reality}, 
the functions $x^i (u_1,\dots,u_n)$, $i=1,\dots,n+1$ are
real--valued. 
Finally, by (8) we obtain that:
$$
ds^2 = \bigl(-C_1 f_1 ^2\bigr)\, du_1 ^2 + \dots + \bigl(- C_n f_n
^2\bigr)\, du_n ^2
$$
is metric of constant curvature 
$K=-1$. So the Lam\'{e} coefficients $H_i =
\sqrt{-C_i f_i ^2}$ satisfy 
the equations (3) and (4).
\qed

%Orthogonal coordinate systems in spaces of %constant curvature corresponding to singular %curves. Examples.

\section{Orthogonal coordinate systems for $S^n$ and $H^n$ corresponding to singular curves. Examples.}

%In this section we show how the modified %Krichever construction
%realizes for the case of a singular spectral %curve $\Gamma$.

In this section we demonstrate the modified Krichever's method introduced in 
\S\,\ref{modified_krichever_sec} 
for the case of a singular 
spectral curve $\Gamma$.
We will consider 
singular reducible curves for which
the irreducible components are 
isomorphic to $\Bbb {CP}^1$ and
the singular points are obtained by
identifying points $a_j$ and $b_j$ that belong 
to different copies of 
$\Bbb {CP}^1$ for $j=1,\dots,s$, 
where $s$ is the number of 
singular points. 
We note that such singular curves appeared in
[7,\,12].

%We will consider the most simple 
%case when a singular curve is reducible and %all irreducible components are
%isomorphic to $\Bbb {CP}^1$. 
%herewith the singular points appear
%as intersections by pairs of points $a_j, %b_j$, $j=1\dots s$ (where $s$ is a number of %singular points) belonging to different 
%copies of $\Bbb {CP}^1$. 

 A Baker--Akhiezer function  is defined %specified 
 on each component 
 $\Bbb {CP}^1$ so that
  its values at 
  the points $a_j,b_j$ coincide.
  We require that the meromorphic differential $\Omega$ in %from 
  Theorem \ref{main_thm} 
  %the theorem in \S\,3 to have 
  has simple poles at 
  the points $a_j$ and $b_j$   
  %and the sum of the corresponding residues to %be equal to zero:
  and the sum of the residues at these 
  points is equal to zero: 
  \begin{equation}
  \label{sum_of_res}
  \operatorname{Res}_{a_j} \Omega + \operatorname{Res}_{b_j} \Omega =
  0, \quad j =1, \dots, s. % \tag 9
  \end{equation}
  % All conditions of the theorem remain %without changes but with the only stipulation
  % that the geometric genus $g$ should be %changed to the arithmetic genus of a %singular surface $\Gamma$. 
  %The claims of the theorem remain valid %because of (9).
  It follows from \eqref{sum_of_res} 
  that  Theorem \ref{main_thm} remains 
  valid for a singular cure (with 
  the geometric genus $g$ changed 
  to the arithmetic genus of the singular 
  curve).

  \begin{picture}(170,115)(-100,-65)
    \qbezier(-60,0)(-60,30)(-10,30) \qbezier(-10,30)(40,30)(40,0)
    \qbezier(40,0)(40,-30)(-10,-30) \qbezier(-10,-30)(-60,-30)(-60,0)
    
    \qbezier(10,0)(10,30)(60,30) \qbezier(60,30)(110,30)(110,0)
    \qbezier(110,0)(110,-30)(60,-30) \qbezier(60,-30)(10,-30)(10,0)
    
    \qbezier(80,0)(80,30)(130,30) \qbezier(130,30)(180,30)(180,0)
    \qbezier(180,0)(180,-30)(130,-30) \qbezier(130,-30)(80,-30)(80,0)
    
    \put(-60,0){\circle*{3}} \put(40,0){\circle*{3}}
    \put(40,0){\circle*{3}} \put(80,0){\circle*{3}}
    \put(10,0){\circle*{3}}
    
    \put(110,0){\circle*{3}} \put(25,24){\circle*{3}}
    \put(25,-24){\circle*{3}} \put(180,0){\circle*{3}}
    \put(95,24){\circle*{3}} \put(95,-24){\circle*{3}}

    \put(-10,33){\shortstack{$\Gamma_1$}}
    \put(55,33){\shortstack{$\Gamma_2$}}
    \put(125,33){\shortstack{$\Gamma_3$}}
    \put(-75,0){\shortstack{$P_1$}} \put(115,0){\shortstack{$Q_1$}}
    \put(45,0){\shortstack{$r$}} \put(-5,0){\shortstack{$P_2$}}
    \put(184,-2){\shortstack{$Q_3$}} \put(65,0){\shortstack{$Q_2$}}
    
    \put(15,20){\shortstack{\small{$a$}}}
    \put(6,-26){\shortstack{\small{$-a$}}}
    \put(35,18){\shortstack{\small{$b$}}}
    \put(30,-26){\shortstack{\small{$-b$}}}
    
    \put(85,20){\shortstack{\small{$c$}}}
    \put(76,-26){\shortstack{\small{$-c$}}}
    \put(105,18){\shortstack{\small{$d$}}}
    \put(100,-26){\shortstack{\small{$-d$}}}
    
    \put(56,-37){\shortstack{\small{$\gamma_1$}}}
    \put(60,-30){\circle*{3}}
    
    \put(130,-37){\shortstack{\small{$\gamma_2$}}}
    \put(134,-30){\circle*{3}}
    
    \put(50,-50){\shortstack{Figure 1}}
    \end{picture} 
  
  \smallskip 
  
  %We now give some examples of construction of %orthogonal coordinate systems in
  %$S^2$ and $H^2$.
  
  Below we show some examples of  
  orthogonal coordinate
  systems on $S^2$ and $H^2$
  obtained by the modified Krichever's method. 
  Let us consider a singular curve $\Gamma$
  that consists of three copies %components 
  of $\Bbb {CP}^1$ intersecting at
  four points. 
%  (an intersection means that 
%  two points of different copies of 
%  $\Bbb {CP}^1$ are identified) 
  In Fig.~1 the points $a,-a$ belong to $\Gamma_1$, the points
  $b,c,-b,-c$ belong to $\Gamma_2$ and the point $d, -d$ belong
  $\Gamma_3$. We put:
  $$
  P_1=\infty\in\Gamma_1,\quad P_2=\infty\in\Gamma_2,\quad
  r=0\in\Gamma_1,
  $$
  $$
  Q_1=0\in\Gamma_2,\quad Q_2=\infty\in\Gamma_3,\quad Q_3=0\in\Gamma_3.
  $$
  %Suppose that 
  The singular curve $\Gamma$ admits a holomorphic involution $\sigma:
  \Gamma\rightarrow\Gamma$, $\sigma(z_i)=-z_i$
  for $z_i \in  \Gamma_i, i=1,2,3$ and an
  antiholomorphic involution  $\tau:\Gamma\rightarrow\Gamma$,
  $\tau(z_i)=\bar{z}_i$ for 
  $z_i \in \Gamma_i$, $i=1,2,3$.     
  %  The points $\gamma_1$ and $\gamma_2$
  %  are real and belong to the components %$\Gamma_2$ and $\Gamma_3$
  %  respectively.
  For $\gamma_1$ and $\gamma_2$ we choose 
  some real points on $\Gamma_2$ and 
  $\Gamma_3$, respectively.

%  The meromorphic 1-form $\Omega$ is defined %by the forms $\omega_i$, 
%  $i=1 \dots 3$ on each component $\Gamma_i$ %correspondingly:
  Now the meromorphic $1$--form $\Omega$
  is given by the forms $\omega_i$, 
  $i=1,2,3$ defined on each 
  component $\Gamma_i$ as follows:   
  $$
  \omega_1=\frac{dz_1}{z_1 \left(z_1^2-a^2\right)}, \quad
  \omega_2=\frac{\bigl(z_2 ^2-\gamma _1 ^2\bigr)\,dz_2} {z_2\bigl(z_2
  ^2-b^2\bigr) \bigl(z_2 ^2-c^2\bigr)}, \quad \omega_3=\frac{\bigl(z_3
  ^2- \gamma _2 ^2\bigr)\,dz_3} {z_3 \bigl(z_3 ^2-d^2\bigr)}.
  $$
%The conditions for $\Omega$ to be regular 
%are as follows:  
  The identities \eqref{sum_of_res} take 
  the form: 
  \begin{equation}
  \gathered
  \operatorname{Res}_a\omega_1+\operatorname{Res}_b\omega_2=0,
  \quad\operatorname{Res}_c\omega_2+\operatorname{Res}_d\omega_3=0,
  \\
  \operatorname{Res}_{-a}\omega_1+\operatorname{Res}_{-b}\omega_2=0,
  \quad\operatorname{Res}_{-c}\omega_2+\operatorname{Res}_{-d}\omega_3=0.
  \endgathered
  %\tag 10
  \end{equation}
  Also, we will impose the additional conditions for the residues of $\Omega$ at
  the points $Q_i$, $i=1,2,3$:
  \begin{equation}
  \operatorname{Res}_{Q_1}\Omega=\operatorname{Res}_{Q_2}\Omega
  =\operatorname{Res}_{Q_3}\Omega. %\tag 11
  \end{equation}
  %Consider the particular solution of (10) and %(11):
  Now let us consider the following solution 
  to the system (10), (11):
  $$
  a=\frac{i}{\sqrt{3}},\quad b=\frac{i}{\sqrt{3}}, \quad c=i,\quad
  d=i,\quad \gamma_1=\frac{1}{\sqrt{3}},\quad \gamma_2=1.
  $$
  %We write out the Baker--Akhiezer function %$\psi$ on each component
  %$\Gamma_1, \Gamma_2, \Gamma_3$, %corresponding to our spectral data:
  Let  $\psi$ be a Baker--Akhiezer function 
  corresponding to the chosen spectral data. 
  We denote by $\psi_1,\psi_2$ and $\psi_3$ 
  the restrictions of $\psi$ to $\Gamma_1, \Gamma_2$ and $\Gamma_3$, respectively. Then $\psi_1, \psi_2$ and $\psi_3$ take the form:      
  $$
  \psi _1(u,v,z_1) = e^{uz_1}f_1(u,v), \quad \psi _2(u,v,z_2) =
  e^{vz_2} \left(f_2(u,v)+\frac{g_2(u,v)}{z_2-\gamma _1}\right),
  $$
  $$
  \psi _3(u,v,z_3) =
  \left(f_3(u,v)+\frac{g_3(u,v)}{z_3-\gamma_2}\right). 
  $$
  %Choose  
  \smallskip
  Let $h =\psi (r) = \frac{1}{\sqrt{3}}$. Then, 
  $$
  -
  \frac{h^2\operatorname{Res}_r\Omega}{\operatorname{Res}_{Q_1}\Omega}=1.
  $$  
  \smallskip 
  Therefore, $f_1(u,v)=\frac{1}{\sqrt{3}}$.  %and the other functions
  The functions 
  $f_2,g_2,f_3$ and $g_3$ 
  %are defined by the compatibility condition %of the function
  %$\psi$ on  $\Gamma_1$, $\Gamma_2$, %$\Gamma_3$:
  can be obtained from the identities: 
  $$
  \psi_1(u,v,a)=\psi_2(u,v,b),\psi_1(u,v,z,-a)=\psi_2(u,v,-b),
  $$
  $$
  \psi_2(u,v,c)=\psi_3(u,v,d), \psi_2(u,v,-c)=\psi_3(u,v,-d).
  $$  
% Consider the following functions of two %variables:
  \setcounter{figure}{1}
  \begin{figure}
  	\begin{center} 
  		\includegraphics[scale=0.56]{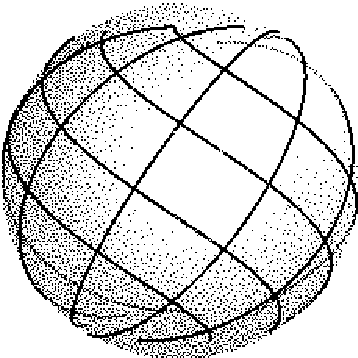}
  	\end{center}
  	\caption{The coordinate lines $u=const$ 
  		     and $v=const$ on $S^2$.}
  	\label{s2_coord_lines_fig}
  \end{figure}
  Finally it follows from Corollary \ref{corollary_sphere} that  
  the functions 
  $x_i (u,v) = \psi_i (u,v,Q_i), i=1,2,3$ 
  provide an orthogonal coordinate system   
  on $S^2$ (see Fig.~\ref{s2_coord_lines_fig}):  
  $$
  x_1(u,v)= \frac{1}{\sqrt{3}} \bigg[\cos\frac{u-v}{\sqrt{3}}+ \sin
  \frac{u-v}{\sqrt{3}}\bigg],
  $$
  \begin{equation*}
  \begin{split}
   x_2(u,v)=\frac{1}{12}\bigg[(3+\sqrt{3})
  \cos\bigg(\frac{u-v}{\sqrt{3}}+v\bigg) +3(-1+\sqrt{3})
  \cos\bigg(\frac{u-v}{\sqrt{3}}-v\bigg)
  \\
  -3(1+\sqrt{3})\sin\bigg(\frac{u-v}{\sqrt{3}}+v\bigg)+(-3+\sqrt{3})
  \sin\bigg(\frac{u-v}{\sqrt{3}}-v\bigg)\bigg],
  \end{split}
  \end{equation*}
  \begin{equation*}
  \begin{split}
  x_3(u,v)=\frac{1}{12}\bigg[3(1+\sqrt{3})
  \cos\bigg(\frac{u-v}{\sqrt{3}}+v\bigg) +(-3+\sqrt{3})
  \cos\bigg(\frac{u-v}{\sqrt{3}}-v\bigg)
  \\
  +(3+\sqrt{3})\sin\bigg(\frac{u-v}{\sqrt{3}}+v\bigg)
  -3(-1+\sqrt{3})\sin\bigg(\frac{u-v}{\sqrt{3}}-v\bigg)\bigg].
  \end{split}
  \end{equation*}
%  Note that if $v$ is fixed then the %coordinate lines are circles with radius
%  $1$, but in contrast to the usual spherical %coordinate system these
%  circles have no common points at the north %and south poles.
  We note that the coordinate lines 
  $v = const$ are circles of radius $1$. 
  These circles do not intersect
  at the north and south poles as 
  for the spherical coordinate system.     
%  By Corollary 1, $u$ and $v$ yield the %orthogonal coordinate system on the %sphere~$S^2$. 
%  Note that in the coordinates $u$, $v$ the %metric on $S^2$ has the form:
  For the coordinates $u$ and $v$ the 
  metric on $S^2$ has the form: 
  $$
  ds^2 = \frac{1}{3}du^2+\frac{1}{3} \left(1-\sin \frac{2 (u-v)}
  {\sqrt{3}}\right)dv^2.
  $$
%  In order to construct an example of an %orthogonal coordinate system in
%  $H^2$ for the same spectral curve $\Gamma$ %(fig. 1), 
We now give an example of an orthogonal 
coordinate system on $H^2$.   
Let us take the
following values for $a,b,c,d,\gamma_1$ and $\gamma_2$: 
  $$
  a = -1,\quad b = 1,\quad c= d = i,\quad \gamma_1 = \sqrt{3},\quad
  \gamma_2 = 1.
  $$
  %It is easy to verify 
  It can be verified 
  that the identities (10) are satisfied. 
  The residues of the form $\Omega$ at the points $Q_1,Q_2,Q_3$ and $r$ are as follows: 
  \begin{equation*}
  	\operatorname{Res}_{Q_1} \Omega = 3, \quad
  	\operatorname{Res}_{Q_2} \Omega = -1,\quad 
  	\operatorname{Res}_{Q_3} \Omega = -1,\quad 
    \operatorname{Res}_{r} \Omega = -1.  
  \end{equation*}	
%  equal to
%  $3$, $-1$, $-1$ and $-1$ respectively.
%   It is convenient for us to put  $h =\psi(r) %= 1$. 
%the compatibility conditions for $\psi$ on the %components $\Gamma_1$,
%$\Gamma_2$, $\Gamma_3$  
  Let $h =\psi(r) = 1$.  
  Therefore,  $f_1(u,v)=1$.  
  The functions $f_2,g_2,f_3$ and $g_3$ 
  can be obtained from the 
  identities:
  $$
  \psi_1(u,v,a)=\psi_2(u,v,b),\psi_1(u,v,z,-a)=\psi_2(u,v,-b),
  $$
  $$
  \psi_2(u,v,c)=\psi_3(u,v,d), \psi_2(u,v,-c)=\psi_3(u,v,-d).
  $$
  \begin{figure}
  \begin{center} 
  	\includegraphics[scale=0.56]{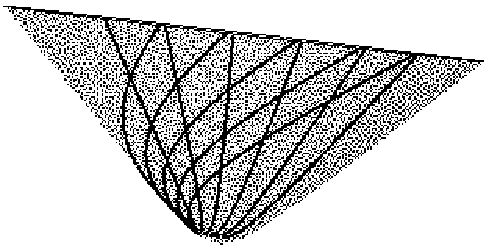}
  \end{center}
    \caption{The coordinate lines $u=const$ 
  		and $v=const$ on $H^2$.}
  	\label{h2_coord_lines_fig}
  \end{figure}
  It follows from Corollary \ref{corollary_hyp_space}
  that  the functions $x_1(u,v) =
  \sqrt{3}\psi(u,v, Q_1)$, 
  $x_2 ( u,v ) = \psi ( u,v,Q_2)$ and 
  $x_3 ( u,v ) = \psi ( u,v,Q_3)$ 
  provide an orthogonal 
  coordinate system on $H^2$ (see Fig.~\ref{h2_coord_lines_fig}):
  $$
  x_1(u,v)= \frac{1}{1+\sqrt{3}} \left[ e^{-u-v}+ (2+\sqrt{3} ) e^{u+v}\right],
  $$
  \begin{equation*}
  \begin{split} x_2(u,v)=\frac{1}{4 \left(1+\sqrt{3}\right)} \bigg[\left(-2e^{-u-v}
  +(6+4 \sqrt{3} ) e^{u+v}\right) \cos v
  \\
  -2\left(\sqrt{3}e^{-u-v} + (2+\sqrt{3}) e^{ u+v}\right)\sin v\bigg],
  \end{split}
  \end{equation*}
  \begin{equation*}
  \begin{split} x_3(u,v)= \frac{1}{4 (2+\sqrt{3})}\bigg[\left(( 3+\sqrt{3} )
  e^{-u-v} + (5+3 \sqrt{3}) e^{u+v}\right) \cos v
  \\
  + \left((-1-\sqrt{3} )e^{-u-v}+(9+5 \sqrt{3} ) e^{u+v} \right) \sin v \bigg].
  \end{split}
  \end{equation*}
%  By Corollary 2, $u$ and $v$ yield the %orthogonal coordinate system
%  in~$H^2$. 
%  Note that in the coordinate system $u$, $v$,
%  the metric on $H^2$ has the form 
  For the coordinates $u$ and $v$ the metric
  on $H^2$ has the form:  
%  $ ds^2 = du^2+( ((7+4 \sqrt{3})e^{-2
%  (u+v)}+ (97+56 \sqrt{3}) e^{2(u+v)})/(52+30 %\sqrt{3})-1 )dv^2. $
  $$
   ds^2 = du^2+\left( 
   \frac{(7+4 \sqrt{3})e^{-2
   	(u+v)}+ (97+56 \sqrt{3}) e^{2(u+v)}}{52+30 \sqrt{3}}-1 \right)dv^2.
  $$ 
 \section*{Acknowledgments}
 
 The authors thank A.~E.~Mironov for helpful conversations and valuable
 comments.

\end{document}